\documentclass[11pt,reqno]{amsart}
\usepackage{amsthm,amssymb,amsbsy}
\usepackage{hyperref}

\input{amssym.def}
\input{amssym}

\makeatletter
\@addtoreset{equation}{section}
\makeatother

\def\be{\begin{equation}}
\def\ee{\end{equation}}
\def\ba{\begin{eqnarray}}
\def\ea{\end{eqnarray}}
\def\lb{\label}

\def\Tr#1{{\rm Tr\str{-1.0}}_{\!R^{\,\mbox{\scriptsize$#1$}}}}
\def\str#1{\rule[#1mm]{0pt}{1mm}}

\def\cal{\mathcal}

\def\tbl#1#2{{\ifmmode \left\{\!\!
\begin{array}{c}
\scriptstyle #1\\[-2pt] \raisebox{2pt}
{$\scriptstyle #2$}
\end{array}
\!\!\right\}
\else
\raisebox{2pt}
{\scriptsize$\left\{\!\!\!
\begin{array}{c}#1\\[-1pt] #2
\end{array}
\!\!\!\right\}$}\fi}}

\def\schf#1{s_{\raisebox{-0.08mm}{$_{#1}$}}}

\newtheorem{theorem}{Theorem}

\newtheorem{lemma}[theorem]{Lemma}
\newtheorem{proposition}[theorem]{Proposition}

\begin{document}
\begin{titlepage}
\title[Power sums of quantum supermatrix]
{Spectral parameterization for the power sums of quantum supermatrix\vspace{1.5cm}}

\author{Dimitri Gurevich}
\address{Dimitri Gurevich, Max Planck Institute for Mathematics,
Vivatsgasse 7, D-53111 Bonn, Germany \&
USTV, Universit\'e de Valenciennes,
59304 Valenciennes, France}
\email{gurevich@univ-valenciennes.fr}

\author{Pavel Pyatov}
\address{Pavel Pyatov,
Max Planck Institute for Mathematics,
Vivatsgasse 7, D-53111 Bonn, Germany \&
Bogoliubov Laboratory of Theoretical Physics, Joint Institute
for Nuclear Research, 141980 Dubna, Moscow Region, Russia}
\email{pyatov@theor.jinr.ru}

\author{Pavel Saponov\vspace{1.5cm}}
\address{Pavel Saponov, Division of Theoretical Physics, IHEP, 142281
Protvino, Moscow region, Russia}
\email{Pavel.Saponov@ihep.ru}
\thanks{The work of D.G. is partially supported by the grant
ANR-05-BLAN-0029-01.
PP and PS are supported by the
Russian Foundation for Basic Research,  grant \mbox{No.\,08-01-00392-a}.
D.G. and P.P. would like acknowledge a warm hospitality of the Max-Planck-Institute f\:ur Mathematik
where this work was completed.
}


\begin{abstract}
A parameterization for the power sums of $GL(m|n)$ type
quantum (super)matrix is obtained in terms of it's spectral values.
\end{abstract}

\maketitle

\end{titlepage}

\section{Introduction}
This paper is a complement to our previous works
\cite{GPS1,GPS2,GPS3} devoted to the quantum matrix
algebras (QMA) of $GL(m|n)$ type.  Here we continue
investigation of the commutative characteristic subalgebra
of the QMA. To be more precise, for the set of `quantum'
traces of `powers' of the quantum (super)matrix we find out
a parameterization in terms of spectral values which are the
quantum analogs of the set of (super)matrix eigenvalues.
Note, that the abovementioned set
of quantum traces generates the characteristic subalgebra.
To illustrate the statements let us briefly recall the corresponding
facts from the classical matrix algebra.

As is well known, any $N\times N$ complex matrix
$M\in {\rm Mat}_N({\Bbb C})$ satisfies a polynomial
Cayley-Hamilton (or, characteristic)  identity,
which can be presented in a factorized  form
$$
\prod_{i=1}^N(M-\mu_i I) = 0\,,
$$
where $I$ is a unit matrix, and
$\mu_i,\,$ $i=1,\dots , N,\,$ are the eigenvalues of $M$.
Opening the brackets one can rewrite the identity in a form
$$
\sum_{k=0}^N(-1)^ke_k(\mu)M^{N-k} = 0\,,
$$
where $e_k(\mu),\,$ $k=0,1,\dots ,N,\,$ are elementary symmetric
polynomials\footnote{For a review of basic results on symmetric
functions see \cite{Mac}.} in  variables $\{\mu_i\}_{1\le i\le N}$.
The elementary symmetric polynomials generate the whole algebra
of symmetric polynomials in the eigenvalues $\mu_i$.
Another well-known generating set for symmetric polynomials
is given by the power sums
$$
p_k(\mu):=\sum_{i=1}^N\mu_i^k\, \equiv\, {\rm Tr}(M^k)\, .
$$
A relation between  two sets of generators is provided by the
Newton's recurrence
$$
k\,e_k+ \sum_{r=1}^{k}(-1)^rp_re_{k-r} = 0\quad \forall\,k\ge 1.
$$

In papers \cite{GPS1,GPS2} analogues of the above classical
results were found for a family of Hecke type QMAs, which includes
the $q$-generalizations of the $GL(m|n)$ type supermatrices
for all integer $m\geq 0$ and $n\geq 0$ (see the next section for
definitions). In particular, the Cayley-Hamilton identities for these
algebras were derived, what clarified the way of
introducing the spectral values for quantum matrices.
It is remarkable that the coefficients of the Cayley-Hamilton polynomials
commute among themselves. They generate a commutative characteristic subalgebra in
the QMA, which serves as an analogue of  the algebra of symmetric polynomials
in the eigenvalues of matrix $M$. It is the Cayley-Hamilton identity which allows one
to present the elements of the characteristic subalgebra as (super)symmetric
polynomials in the spectral values of the quantum matrix.

As in the classical case, the quantum analogs of the power sums
can be defined as some specific traces of `powers' of the quantum
matrix. By their construction the power sums belong to the characteristic subalgebra,
but their explicit expression  in terms of the spectral values was not known yet.
The main goal of the present paper is to derive such an expression.

Our presentation is strongly based on the previous works cited
above. In the next section we give a list of notations, definitions
and main results which will be used below. For more detailed
exposition, proofs and a short overview the reader is referred to \cite{GPS1,GPS2}.

\section{Some basic results and definitions}

Let $V$ be a finite dimensional linear space over the field of
complex numbers $\Bbb C$, $\dim V = N$. Let $I$
denote the identity matrix (its dimension being clear from the
context, if not explicitly specified), and $P\in {\rm Aut}(V^{\otimes 2})$  be the
permutation automorphism: $P(u\otimes v) = v\otimes u$.

With any element
$X\in {\rm End}(V^{\otimes p})$, $p= 1,2,\dots ,$ we associate
a sequence of endomorphisms $X_i\in {\rm End}(V^{\otimes k})$,
$\,k\geq p$, $\,i=1,\dots ,k-p+1$,\, according to the rule
\be
\nonumber
X_i = I^{\otimes {(i-1)}}\otimes X\otimes I^{\otimes ^{(k-p-i+1)}}\,,
\quad 1\leq i\leq k-p+1,
\ee
where $I$ is the identical automorphism of $V$.

Consider a pair of invertible operators $R,F\in {\rm Aut}
(V^{\otimes 2})$ subject to the following conditions:
\begin{enumerate}
\item\label{i-1}
The operators  $R$ and $F$ satisfy the {\em Yang-Baxter equations}
\be
\lb{YB}
R_1R_2R_1 = R_2R_1R_2\, ,\qquad
F_1F_2F_1 = F_2F_1F_2\, .
\ee
Such operators are called {\it {\rm R}-matrices}.
\item\label{i-2}
The pair of R-matrices $\{R,F\}$ is {\em compatible}, that is
\begin{equation}
R_1F_2F_1 = F_2F_1R_2\, ,
\qquad
F_1F_2R_1 = R_2F_1F_2\, .
\label{YBE}
\end{equation}
\item\label{i-3}
The matrices of both operators $R$ and $F$ are {\em strictly
skew invertible}. Taking the operator $R$ as an example, this
requirement means the following:\smallskip
\item[a)]
$R$ is skew invertible if there exists an operator $\Psi^R\in {\rm End}(V^{\otimes 2})$
such that
$$
{\rm Tr}_{(2)}R_{12}\Psi_{23}^R  = P_{13}\, ,
$$
where the subscript in the notation of the trace shows the
number of the space  $V$,  where the trace is evaluated
(the enumeration of the component spaces in the tensor product is taken as follows
$V^{\otimes k} := V_1\otimes V_2\otimes\dots  \otimes V_k$).
\item[b)]
The strictness condition implies additionally that the operator
$D_1^R:={\rm Tr}_{(2)}\Psi^R_{12}$ is invertible.
\end{enumerate}
With the matrix $D^R$ one  defines the {\em R-trace}
operation\footnote{In a literature on quantum groups the
R-trace is usually named the quantum trace or, shortly, the
$q$-trace. Giving the different name to this operation we hope
to avoid misleading associations with a parameter $q$ of the
Hecke algebra (see below).}
$\Tr{}:{\rm Mat}_N(W)\rightarrow W$
$$
\Tr{}(X) := \sum_{i,j=1}^N{D^R}_i^jX_j^i,\quad  X\in
{\rm Mat}_N(W),
$$
where $W$ is any linear space.
\smallskip

Given a compatible pair $\{R,F\}$ of strictly skew invertible
R-matrices
the {\em quantum matrix algebra } ${\cal M}(R,F)$
is a unital associative  algebra  generated
by $N^2$ components of the matrix  $\|M^i_j\|_{i=1}^N$ subject
to the relations
\be
\label{def:QM}
R_1M_{\overline 1}M_{\overline 2} = M_{\overline 1}
M_{\overline 2} R_1\, .
\ee
Here we have introduced a notation
\be
\nonumber
M_{\overline 1} := M_1, \quad M_{\overline {k+1}}  :=
F_k M_{\overline k} F^{-1}_k
\ee
for the copies $M_{\overline{k}}$ of the matrix $M$. The defining
relations (\ref{def:QM}) then imply the same type relations for
consecutive pairs of the copies of $M$ (see  \cite{IOP2})
\be
R_k\,M_{\overline k}M_{\overline{k+1}} =
M_{\overline k}M_{\overline{k+1}}\, R_k.
\nonumber
\ee

Specific subfamilies in a variety of QMAs are extracted by imposing
additional conditions on the R-matrix $R$ in the
definition (\ref{def:QM}). These conditions we are now going to describe.
\smallskip

Suppose that $R$ has a quadratic minimal polynomial
which can be suitably normalized as\footnote{Note that all the conditions
(\ref{YB}), (\ref{YBE}) do not depend on
normalization of $R$.}
\be
(R-q\, I)(R+q^{-1}\,I)=0\,, \quad q\in {\Bbb C}\setminus 0\,.
\label{Hec}
\ee
This relation in the present context is called the {\em Hecke condition}
and the R-matrices satisfying it are called {\it the
Hecke }  R-matrices. We further assume that the  parameter $q$ in (\ref{Hec}) is
{\em generic}, that means, it does not coincide with the roots
of equations
\begin{equation}
k_q:={q^k-q^{-k}\over q-q^{-1}}= 0\qquad \forall\,k=2,3,\dots\,.
\label{q-gen}
\end{equation}

Given any Hecke R-matrix $R$, one can construct a series of
{\em {\rm R}-matrix representations} $\rho_R$ of the A type Hecke algebras
${\cal H}_k(q)\stackrel{\rho_R}{\longrightarrow} {\rm End} (V^{\otimes k})$,
$k=2,3,\dots$. The characteristic properties of these representations
are used for a classification of the Hecke R-matrices\,\footnote{
A brief description of the Hecke algebras and their R-matrix
representations can be found in \cite{GPS1}. For a more detailed exposition of
the subject the reader is referred to \cite{R,OP1} and to the references therein.}.
Not going into details of the construction we only mention that under conditions (\ref{q-gen})
the Hecke algebra ${\cal H}_k(q)$ is isomorphic to the group algebra of the symmetric group
${\Bbb C}[S_k]$ and its irreducible representations are labelled by a set of
partitions $\lambda\vdash k$, the corresponding central idempotents in ${\cal H}_k(q)$ are
further denoted as $e^\lambda$. We fix some decomposition of $e^\lambda$
into the sum of primitive idempotents
$e^\lambda_a\in{\cal H}_k(q):\; e^\lambda =\sum_{a=1}^{d_\lambda}e^\lambda_a$,
where $d_\lambda$ is the dimension of the representation with label $\lambda$.
It is also suitable to introduce the following notations:
\begin{enumerate}
\item[--\,]
Given two arbitrary integers $m\ge 0$ and $n\ge 0$, an infinite
set of partitions $\lambda =
(\lambda_1,\lambda_2,\dots)$, satisfying restriction
$\lambda_{m+1}\le n$ is denoted as $\mbox{\sf H}(m,n)$.
\item[--\,]
The partition $((n+1)^{m+1})\vdash (m+1)(n+1)$
is shortly denoted as $\lambda_{m,n}$.
The corresponding Young diagram is a rectangle with
$m+1$ rows of the length $n+1$. Note that
$\lambda_{m,n}$ is a minimal partition not belonging to the set
$\mbox{\sf H}(m,n)$.
\end{enumerate}

Now we are ready to formulate the classification of the Hecke R-matrices.

\begin{proposition}{\bf(\cite{H,GPS3})}
\lb{prop-1}
For a generic value of $q$ the set of the Hecke {\rm R}-matrices is
separated into subsets
labelled by an ordered pair of non-negative integers $\{m,n\}$. {\rm R}-matrices
belonging to the subset with label $\{m,n\}$ are called {\em $GL(m|n)$ type} ones
{\rm (}alternatively, they are assigned a {\em bi-rank} $(m|n)${\rm )}.
{\rm R}-matrix representations $\rho_R$  generated by a  $GL(m|n)$ type
{\rm R}-matrix $R$ fulfill the following criterion:
for all integer $k\geq 2\,$ and for any partition $\nu\vdash k$
the images of the idempotents $e^{\nu}\in {\cal H}_k(q)$ satisfy the relations
$$
\rho_R(e^\nu) = 0 \quad\mbox{iff} \quad
\nu\not\in \mbox{\sf H}(m,n),
$$
or, equivalently, iff $\lambda_{m,n}\subset \nu$,
where the inclusion $\mu=(\mu_1,\mu_2,\dots)\subset \nu =
(\nu_1,\nu_2,\dots)$ means that $\mu_i\le \nu_i$ $\forall\, i$.
\end{proposition}

The algebra ${\cal M}(R,F)$ defined by a Hecke ($GL(m|n)$ type) R-matrix $R$
is further referred to as the {\em Hecke ($GL(m|n)$ type)} quantum matrix algebra.
\medskip

For the Hecke type QMA  ${\cal M}(R,F)$ we consider
a set of its elements $s_\lambda(M)$ called
{\em the Schur functions}
\begin{equation}
s_0(M):= 1,\qquad s_\lambda(M):=\Tr{(1\dots p)}
(M_{\overline 1}\dots M_{\overline p}\,
\rho_R(e^\lambda_a))\, ,
\quad \lambda\vdash p,\quad p=1,2,\dots\, ,
\nonumber
\end{equation}
where the latter formula does not depend on a particular choice of the primitive idempotent
$e^\lambda_a$ (actually, one can substitute it by $d_\lambda^{-1}e^\lambda$).
As was shown in \cite{IOP1}, a linear
span of the Schur functions $s_\lambda(M)\;\forall\,\lambda$,
is an abelian subalgebra in ${\cal M}(R,F)$. We further call it the
{\em characteristic subalgebra} of ${\cal M}(R,F)$.
It  follows that
the  characteristic subalgebra  of the
$GL(m|n)$ type QMA is spanned  by the
Schur functions $s_\lambda(M)$, $\lambda\in \mbox{\sf H}(m,n)$.
The multiplication table for the elements $s_\lambda(M)\in {\cal M}(R,F)$
coincides with that for the basis of Schur functions in the ring of symmetric
functions (see \cite{Mac}) thus justifying the notation. One has \cite{GPS2}
\be
\lb{LR}
s_\lambda(M) s_\mu(M) = \sum_{\nu} C_{\lambda\mu}^\nu s_\nu(M),
\ee
where $C_{\lambda\mu}^\nu$ are the Littlewood-Richardson coefficients.
Later on we shall need  the information about generating sets of the
characteristic subalgebra.

\begin{proposition}
{\bf (\cite{IOP1,IOP2})} For generic values of $q$ the characteristic subalgebra
of the Hecke QMA  ${\cal M}(R,F)$ is generated by any one
of the following three sets
\begin{enumerate}
\item
the single column Schur functions:
$\,a_k(M):= s_{(1^k)}(M)$, $\,k= 0,1,2,\dots$;\vspace{1mm}
\item
the single row Schur functions:
$\,s_k(M):= s_{(k)}(M)$, $\,k= 0,1,2,\dots$;\vspace{1mm}
\item
the set of {\em power sums}:
\begin{equation}
p_{\,0}(M):=(\Tr{}I)\,1,\quad
p_k(M):= \Tr{(1\dots k)}(M_{\overline 1}\dots
M_{\overline k}\,R_{k-1}\dots R_1)\,, \quad k\ge 1\,.
\lb{st-sm}
\end{equation}
\end{enumerate}
These sets are connected by a series of recursive {\em Newton and Wronski relations}
\begin{eqnarray}
(-1)^k k_q\, a_k(M) + {\textstyle \sum_{r=0}^{k-1}}\, (-q)^{r}a_r(M)\,
p_{k-r}(M) &=& 0\, ,\label{q-anti}\\[1mm]
k_q\, s_k(M) - {\textstyle \sum_{r=0}^{k-1}}\, q^{-r}s_r(M)\,
p_{k-r}(M) &=& 0\,, \label{q-simm}\\[1mm]
{\textstyle \sum_{r=0}^k}\, (-1)^r a_r(M)\, s_{k-r}(M) &=& 0 \qquad \forall\,k\ge 1\,.
\label{wronski}
\end{eqnarray}
\end{proposition}\smallskip

On introducing the generating functions for these sets of generators
\be
\lb{gen-func}
A(t):=\sum_{k\geq 0} a_k(M)\, t^k ,\quad
S(t):=\sum_{k\geq 0} s_k(M)\, t^k ,\quad
P(t):=1+(q-q^{-1})\sum_{k\ge 1}p_k(M)\,t^k\, ,
\ee
one can rewrite relations (\ref{q-anti})--(\ref{wronski}) in a compact form
\cite{Mac,I}\footnote{The first of these relations in the limiting case
${\cal M}(R,R)\stackrel{q\rightarrow 1}{\longrightarrow} U({\mathfrak g\mathfrak l}_n)$
was also derived in \cite{U}.}
\be
P(-t)A(qt) = A(q^{-1}t)\,,
\qquad
P(t)S(q^{-1}t) = S(qt)\,,
\qquad A(t)S(-t)=1\, .
\label{newton-2}
\ee
For the $GL(m|n)$ type QMA the zeroth power sum equals \cite{GPS3}
\begin{equation}
p_{\,0}(M) = q^{n-m}(m-n)_q 1\,.
\label{tr-id}
\end{equation}

\vspace{2mm}
One of the remarkable properties of the Hecke QMA ${\cal M}(R,F)$ is the existence of the
characteristic identity for the matrix $M$ of it's generators.
To formulate the result we introduce a notion of the {\em matrix $\star$-product}  of
{\em quantum matrices} (for a detailed exposition see \cite{OP2}, section 4.4).
Namely, starting with the quantum matrix of generators $M$
and the scalar quantum matrices $s_\lambda(M) I\;\; \forall\,\lambda\vdash k,\; k\geq 0,\,$
we construct the whole set of quantum matrices by the following recursive procedure:
given any quantum matrix $N$, its $\star$-multiplication by $s_\lambda(M) I$ and left
$\star$-multiplication by $M$ are also quantum matrices defined as
\ba
\nonumber
&&M\star (s_\lambda(M) I)=(s_\lambda(M) I)\star M := M\cdot s_\lambda(M) ,
\\
\nonumber
&&M\star N := M\cdot \phi(N), \quad \mbox{where}\quad \phi(N)_1:=\Tr{(2)}N_{\overline{2}}R_{12},
\ea
Here the dot-product means the usual multiplication of a matrix by a scalar(matrix).
The $\star$-product of the quantum matrices is demanded to be associative and,
by construction, it is commutative. The $\star$-powers of the quantum matrix $M$
read
\begin{equation}
M^{\overline 0}:= I, \quad
M^{\overline 1}:= M,\quad M^{\overline k}:=
\raisebox{-5mm}{$\stackrel{\underbrace{M\star\dots\star M}}{\mbox{\tiny$k$ times}}$} =
\Tr{(2\dots k)}
(M_{\overline 1}\dots M_{\overline k}
\,R_{k-1} \dots R_1)\;\forall\,  k>1.
\nonumber
\end{equation}
We note that for the family of the so-called {\it reflection equation
algebras}  --- these are the QMAs of the form ${\cal M}(R,R)$ --- the
$\star$-product is identical to the usual matrix product.
\medskip

The characteristic identity depends essentially on a type of the quantum matrix algebra.
For the $GL(m|n)$ type QMA  it is an $(m+n)$-th order polynomial identity
in $\star$-powers of the matrix $M$ with coefficients in the characteristic subalgebra.
One has the following $q$-analogue of the classical Cayley-Hamilton theorem.

\begin{theorem}
{\bf (\cite{GPS1,GPS2})}
The characteristic identity for
the matrix of generators of the $GL(m|n)$ type QMA ${\cal M}(R,F)$ reads
\be
\Big(\sum_{k=0}^m (-q)^k\, M^{\overline{m-k}}
\schf{[m|n]^k}(M)
\Big) \star \Big(\sum_{r=0}^n q^{-r}\, M^{\overline{n-r}}
\schf{[m|n]_r}(M)\Big)\equiv 0\, ,
\label{factor-ch}
\ee
where we used a shorthand notation for the partitions
$$
[m|n]^k := \bigl( (n+1)^k, n^{m-k}\bigr), \qquad
[m|n]_r := \bigl( n^m, r\bigr).
$$
\end{theorem}

Remarkably enough, for generic type QMA (i.e., if $mn>0$)
the characteristic identity (\ref{factor-ch}) factorizes in
two parts.
Therefore, when setting factorization problem for the characteristic
polynomial one is forced to separate all the $\,m+n\,$ roots into two
parts of sizes $m$ and $n$.

Let ${\Bbb C}[\mu,\nu]$ be an algebra   of
polynomials in two sets of mutually commuting and algebraically
independent variables $\mu:=\{\mu_i\}_{1\le i\le m}$ and
$\nu:=\{\nu_j\}_{1\le j\le n}$.  Consider a map of the
coefficients of the characteristic polynomial into
${\Bbb C}[\mu,\nu]$
\begin{eqnarray}
\schf{[m|n]^k}(M)
&\mapsto &
\schf{[m|n]^k}(\mu,\nu)
:=
\schf{[m|n]}(\mu,\nu)\, e_k(q^{-1}\mu)\,
, \quad 1\le k\le m\, ,
\label{def:mu}
\\[1mm]
\schf{[m|n]_r}(M)
&\mapsto&
\schf{[m|n]_r}(\mu,\nu)  := \schf{[m|n]}(\mu,\nu)\, e_r(-q\nu)\, ,
\quad\;\, 1\le r\le n\, ,
\label{def:nu}
\end{eqnarray}
where $e_k(\cdot)$ are the elementary symmetric polynomials in
their arguments
(e.g., $e_k(\mu)\equiv e_k(\mu_1,\dots ,\mu_m)=
\sum_{1\le i_1<\dots <i_k\le m}\mu_{i_1}\dots \mu_{i_k}$\vspace{1mm}). For the moment
we don't specify an explicit expression for the polynomial $\schf{[m|n]}(\mu,\nu)$.
We now define a central extension of the $\star$-product algebra of the quantum matrices
by the scalar matrices of the form $p(\mu,\nu) I$,  $\,p(\mu,\nu)\in {\Bbb C}[\mu,\nu]$,\,
such that $s_\lambda(M) I\equiv s_\lambda(\mu,\nu) I$.
In the extended algebra the characteristic identity (\ref{factor-ch}) takes a completely
factorized form
\be
\prod_{i=1}^m
(M - \mu_i I) \star
\prod_{j=1}^n
( M - \nu_j I)\cdot (\schf{[m|n]}(\mu,\nu))^{2}
\equiv 0\, .
\nonumber
\ee
Assuming that $s_{[m|n]}(\mu,\nu)\neq 0$ we can interpret
the variables $\mu_i, \;i=1,\dots ,m,\,$ and $\nu_j,\; j=1,\dots ,n,\,$
as eigenvalues of the quantum matrix $M$. They are called, respectively,
{\em ``even''} and {\em ``odd'' spectral values} of $M$.

The map (\ref{def:mu}), (\ref{def:nu}) admits a unique extension
to a homomorphic map of the characteristic subalgebra into the algebra
${\Bbb C}[\mu,\nu]$ of polynomials in spectral values $\mu_i$ and $\nu_j$.
Using the Littlewood-Richardson multiplication rule (\ref{LR}) we obtain
(see \cite{GPS2})
\begin{eqnarray}
a_k(M)\equiv s_{[k|1]}(M)&\mapsto &
a_k(\mu,\nu)\, :=\, \sum_{r=0}^ke_r(q^{-1}\mu)\,h_{k-r}(-q\nu)\, ,
\label{sk-a}\\
s_k(M)\equiv s_{[1|k]}(M)&\mapsto &s_k(\mu,\nu)
\, :=\,  \sum_{r=0}^ke_r(-q\nu)\,h_{k-r}(q^{-1}\mu)\, ,
\label{sk-s}
\end{eqnarray}
where $\,h_k(\dots)\,$ stands for the complete symmetric polynomial
in its variables:
$h_k(\mu)\equiv $ $h_k(\mu_1,\dots ,\mu_m) =
\sum_{1\leq i_1\leq\dots \leq i_k\leq m}\mu_{i_1}\dots \mu_{i_k}$\vspace{1mm}.
Since each of the sets $\{a_k(M)\}_{k\geq 0}$,\, $\{s_k(M)\}_{k\geq 0}$ generates the
characteristic subalgebra, the homomorphism is completely defined by (\ref{sk-a}) or
by (\ref{sk-s}).
In particular, formulas (\ref{sk-a}),  (\ref{sk-s}) prescribe an explicit
expression for the unspecified polynomial $s_{[m|n]}(\mu,\nu)$
in (\ref{def:mu}), (\ref{def:nu}),
which is the image of $\schf{[m|n]}(M)$:
\begin{equation}
\schf{[m|n]}(M)\, \mapsto\, s_{[m|n]}(\mu,\nu) =
\prod_{i=1}^m\prod_{j=1}^n \left(q^{-1}\mu_i - q\nu_j\right).
\nonumber
\end{equation}

This homomorphic map induced by (\ref{sk-a}), or (\ref{sk-s}) is called the
{\em spectral parameterization} of the characteristic subalgebra.
In the next section we derive the spectral parameterization for the third generating set
of the characteristic subalgebra,
that is the set of power sums $\{p_k(M)\}_{k\geq 0}$.

\section{Spectral parameterization of power sums}

In this section we are working with the $GL(m|n)$ type QMA ${\cal M}(R,F)$
defined by relations (\ref{def:QM}) with an R-matrix $R$ satisfying the criterion
of the proposition \ref{prop-1}. During the considerations we
assume that the parameter $q$ is generic (see (\ref{q-gen})), although afterwards
this restriction can be waived out: unlike the cases of $a_k(M)$ and $s_k(M)$
the power sums $p_k(M)$ are consistently defined for all $q\in {\Bbb C}\setminus 0$.

For the particular case of the $GL(m)\equiv GL(m|0)$ type
reflection equation algebra  ${\cal M}(R,R)$
the spectral parameterization
of the power sums was found in \cite{GS}. Taking into account that in the
$GL(m)$ case the quantum matrix has only ``even'' eigenvalues
$\{\mu_i\}_{1\le i\le m}$, the result reads
$$
p_k(M)\mapsto \sum_{i=1}^m d_i\, \mu_i^k\, ,
\qquad\mbox{where}\qquad
d_i :=q^{-1} \prod_{j\not=i}^m{{\mu_i-q^{-2}\mu_j}
\over{\mu_i-\mu_j}}\,.
$$

Our goal is to extend the this formula for the general case.
To this end we introduce an auxiliary set of polynomials
$\{\pi_k(\mu,\nu)\}_{k\ge 1}\subset {\Bbb C}[\mu,\nu]$
defined by  relations
\begin{equation}
\pi_k(\mu,\nu) := \sum_{i=1}^m(q^{-1}\mu_i)^k -
\sum_{j=1}^m(q\nu_j)^k\,,\qquad k\ge 1\,.
\label{def:pi}
\end{equation}
In the subsequent considerations we shall use the following property
of these polynomials.
\begin{lemma}
The sets of polynomials
$\{a_k(\mu,\nu)\}_{k\geq 0}$, $\{s_k(\mu,\nu)\}_{k\geq 0}$ (see (\ref{sk-a}), (\ref{sk-s}))
and $\{\pi_k(\mu,\nu)\}_{k\geq 1}$ satisfy the Newton's recurrent relations
\begin{eqnarray}
(-1)^k k\, a_k(\mu,\nu) + \sum_{r=0}^{k-1} (-1)^{r}
a_r(\mu,\nu)\, \pi_{k-r}(\mu,\nu) &=& 0\, ,
\label{a-pi}\\[-1mm]
k\, s_k(\mu,\nu) - \sum_{r=0}^{k-1}s_r(\mu,\nu)\,
\pi_{k-r}(\mu,\nu) &=& 0\,
\quad \forall\,k\ge 1\,.
\label{s-pi}
\end{eqnarray}
\end{lemma}

\noindent{\bf Proof.} First, we recall  relations among the generating functions
for the power sums and the elementary symmetric and complete symmetric
polynomials in a finite set of variables $x:=\{x_i\}_{1\leq i\leq p}$
(see \cite{Mac}, section I.2)
\ba
\nonumber
E(x|t)&:=&\sum_{k=0}^p e_k(x) t^k\, =\, \prod_{i=1}^p(1+x_it)\, ,\quad H(x|t)\, :=\,
\sum_{k\geq 0} h_k(x) t^k\, =\,
\prod_{i=1}^p(1-x_it)^{-1}\,,\hspace{8mm}
\\
\lb{P-EH}
P(x|t)&:=&\sum_{k\geq 1}p_k(x) t^{k-1}\, =\, - {d\over dt}\log E(x|-t)\, =\,
{d\over dt}\log H(x|t)\, .
\ea
Consider three functions depending on two sets of variables
$x:=\{x_i\}_{1\leq i\leq m}$ and $y:=\{y_i\}_{1\leq i\leq n}$:
\ba
\nonumber
A(x,y|t)&:=& E(x|t) H(-y|t), \qquad S(x,y|t)\,:=\, H(x|t) E(-y|t),
\\
\lb{Pi}
\Pi(x,y|t)&:=& P(x|t) - P(y|t).
\ea
These functions serve as super-matrix analogues of, respectively, generating functions of the
elementary and complete symmetric polynomials and the power sums
(see \cite{Mac}, section I.5, exercise 27 and bibliographic references for it).
Indeed, using (\ref{P-EH}) it is easy to check that these functions
satisfy  relations similar to (\ref{P-EH})
\be
\lb{Pi-AS}
\Pi(x,y|t)\, =\, - {d\over dt}\log A(x,y|-t)\, =\,
{d\over dt}\log S(x,y|t)\, .
\ee
Now the assertion of the lemma follows from an observation that
$A(q^{-1}\mu, q\nu |t)$, $S(q^{-1}\mu, q\nu |t)$ and $\Pi(q^{-1}\mu, q\nu |t)$
are, respectively, generating functions for the sets of polynomials
$\{a_k(\mu,\nu)\}_{k\geq 0}$, $\{s_k(\mu,\nu)\}_{k\geq 0}$ and $\{\pi_k(\mu,\nu)\}_{k\geq 1}$.
The relations (\ref{a-pi}) and (\ref{s-pi}) are just expansions of (\ref{Pi-AS})
in powers of $t$.
\hfill\rule{6.5pt}{6.5pt}
\medskip

Now we can formulate the main result on the power sums.
\begin{proposition}
 The spectral parameterization (\ref{sk-a}) (or (\ref{sk-s})) for the power sums
$p_k(M)$ (\ref{st-sm}) in the $GL(m|n)$ type quantum matrix algebra is given by
the formulas
\begin{equation}
p_k(M)\mapsto p_k(\mu,\nu)=\sum_{i=1}^m d_i \mu_i^k +
\sum_{j=1}^n {\tilde d}_j\nu_j^k\quad \forall\,k\ge 0\,,
\label{pk}
\end{equation}
where  the ``weight'' coefficients $d_i$ and
$\tilde d_j$  explicitly read
\begin{eqnarray}
d_i &:=& q^{-1}\prod_{p=1\atop p\not=i}^m
\frac{\mu_i - q^{-2}\mu_p}{\mu_i-\mu_p}\,
\prod_{j=1}^n \frac{\mu_i - q^2\nu_j}{\mu_i -\nu_j}\, ,
\label{d-i}\\
\tilde d_j &:=& -\,q\,\prod_{i=1}^n \frac{\nu_j - q^{-2}\mu_i}{\nu_j-\mu_i}\,
\prod_{p=1\atop p\not=j}^n \frac{\nu_j - q^2\nu_p}{\nu_j-\nu_p}\,.
\label{dd-i}
\end{eqnarray}
Recall, that the spectral values
$\{\mu_i\}$ and $\{\nu_j\}$ are supposed to be algebraically
independent, and therefore all the coefficients $d_i$ and $\tilde d_j$
are nonzero and well defined.
\end{proposition}

\noindent
{\bf Proof.}
For the proof we
need yet another recursive set of formulas for the power
sums $\{p_k(\mu,\nu)\}_{k\ge 1}$ (\ref{pk}) and the polynomials
$\{\pi_k(\mu,\nu)\}_{k\ge 1}$ (\ref{def:pi}).
\begin{lemma}
The following relations hold true
\begin{equation}
k\, p_k(\mu,\nu) = k_q\,\pi_k(\mu,\nu) + (q-q^{-1})\sum_{r=1}^{k-1}
r_q\,\pi_r(\mu,\nu)p_{k-r}(\mu,\nu)\quad \forall\,k\ge 1\,.
\label{p-pi}
\end{equation}
In terms of the generating functions
$$
P(\mu,\nu|t):=1+(q-q^{-1})\sum_{k\geq 1}p_k(\mu,\nu)t^k\, \qquad
\Pi(q^{-1}\mu, q\nu|t):=\sum_{k\geq 1}\pi_k(\mu,\nu) t^{k-1}
$$
($P(\mu,\nu|t)$ is the spectral parameterization of $P(t)$, see (\ref{gen-func});
$\Pi(\mu,\nu|t)$ was first defined  in (\ref{Pi}))
the relations  (\ref{p-pi}) shortly read
\be
P(\mu,\nu|t)\Bigl(q\,\Pi(q^{-1}\mu,q\nu|qt)\, -\, q^{-1}\,\Pi(q^{-1}\mu,q\nu|q^{-1}t)\Bigr)
= \frac{d}{dt}P(\mu,\nu|t)\,.
\lb{P-Pi}
\ee
\end{lemma}
\noindent
{\bf Proof of the lemma.}
Introduce a meromorphic function $f:{\Bbb C}[\mu,\nu]
\rightarrow {\Bbb C}(\mu,\nu)$ by the formula
$$
f(z):= \prod_{i=1}^m \frac{(z - q^{-2}\mu_i)}{(z-\mu_i)}
\prod_{j=1}^n \frac{(z - q^2\nu_j)}{(z -\nu_j)}\,.
$$
Since the spectral values are algebraically independent the
above function has a first order pole  at
each point $\mu_i$ and $\nu_j$. Besides, as can be easily seen
$$
f(0) = q^{2(n-m)}\qquad\mbox{and}\qquad
\lim_{z\rightarrow\infty}f(z) = 1\,.
$$
Taking into account the above limit at infinity, we expand
the function $f(z)$ into the sum of simple fractions
\begin{equation}
f(z) = 1 + \sum_{i=1}^m\frac{1}{z-\mu_i}\,
{\rm Res}f(z)_{\rule[-1mm]{0.4pt}{4mm}_{\,z=\mu_i}}
+ \sum_{j=1}^n\frac{1}{z-\nu_j}\,
{\rm Res}f(z)_{\rule[-1mm]{0.4pt}{4mm}_{\,z=\nu_j}}\,,
\label{si-frac}
\end{equation}
where the residues at the poles read
$$
{\rm Res}f(z)_{\rule[-1mm]{0.4pt}{4mm}_{\,z=\mu_i}} =
(q-q^{-1})\,\mu_i\,d_i\,,\qquad
{\rm Res}f(z)_{\rule[-1mm]{0.4pt}{4mm}_{\,z=\nu_j}} =
(q-q^{-1})\,\nu_j\,\tilde d_j\,.
$$
Evaluating the right hand side of (\ref{si-frac}) at $z=0$ we get
$$
f(0) = q^{2(n-m)} = 1 - (q-q^{-1})\Bigl(\sum_{i=1}^m d_i+
\sum_{j=1}^n\tilde d_j\Bigr):=1-(q-q^{-1})\,p_{\,0}(\mu,\nu)
$$
and therefore,
$$
p_0(\mu,\nu) = q^{n-m}(m-n)_q 1\,.
$$
So, we have verified the consistency of (\ref{pk}) with our previous result
(\ref{tr-id}) for $p_{\,0}(M)$.

In order to prove relations (\ref{p-pi}) we  expand
functions
$z^kf(z)$, $\,k\ge 1,\,$ into the sum of simple fractions.
Besides the simple poles
in $\mu_i$ and $\nu_j$, the function $z^kf(z)$  possesses the $k$-th
order pole
at $z=\infty$, or, introducing a new variable $y=z^{-1}$, at the point
$y=0$.
Taking into account that
$$
{\rm Res}\,z^kf(z)_{\rule[-1mm]{0.4pt}{4mm}_{\,z=\mu_i}}=
(\mu_i)^k\,{\rm Res}f(z)_{\rule[-1mm]{0.4pt}{4mm}_{\,z=\mu_i}}
\,,\qquad
{\rm Res}\,z^kf(z)_{\rule[-1mm]{0.4pt}{4mm}_{\,z=\nu_j}} =
(\nu_j)^k\,{\rm Res}f(z)_{\rule[-1mm]{0.4pt}{4mm}_{\,z=\nu_j}}
\,,
$$
we come to the corresponding expansion
\be
z^kf(z) = \sum_{r=0}^{k}\frac{z^{k-r}}{r!}\, f_r + (q-q^{-1})\Bigl(
\sum_{i=1}^m\frac{d_i\mu_i^{k+1}}{z-\mu_i} +
\sum_{j=1}^n\frac{\tilde d_i\nu_j^{k+1}}{z-\nu_i}\Bigr)\, ,
\label{zf-exp}
\ee
where the coefficients $f_r$ are the following derivatives
$$
f_r:=\frac{d^r\!f(y)}{dy^r}\,\rule[-2.5mm]{0.4pt}{5mm}_{\,y=0}\,.
$$
Evaluating (\ref{zf-exp}) at the point $z=0$ and taking into account
(\ref{pk}) we find the relation
\begin{equation}
p_k(\mu,\nu) = \frac{f_k}{(q-q^{-1})\,k!}\quad\forall\, k\ge 1\,.
\label{p-f}
\end{equation}
Let us turn to the calculation of $f_k$.
Recalling that $y=z^{-1}$ we
write the first order
derivative $f^\prime(y)$
in the form
\begin{equation}
\frac{df(y)}{dy} = (q-q^{-1})\,u(y)f(y)\,,
\label{f-prime}
\end{equation}
where the function $u(y)$ reads
$$
u(y) := \sum_{i=1}^m \frac{q^{-1}\mu_i}{(1-\mu_iy)(1-q^{-2}\mu_iy)}
-\sum_{j=1}^n \frac{q\nu_j}{(1-\nu_jy)(1-q^2\nu_jy)}\,.
$$
By a simple induction in $k$ one can check that
\begin{eqnarray*}
u_k(y)\,:=\,\frac{d^ku(y)}{dy^k} &=&k!\Bigl\{\sum_{i=1}^m\sum_{a=0}^k
\frac{q^{-1-2(k-a)}\mu_i^{k+1}}{(1-\mu_iy)^{1+a}(1-q^{-2}\mu_iy)^{1+k-a}}\hspace{20mm}
\\[.5mm]
&&\hspace{7mm}-\, \sum_{j=1}^n\sum_{a=0}^k \frac{q^{1+2(k-a)}\nu_j^{k+1}}
{(1-\nu_jy)^{1+a}(1-q^2\nu_jy)^{1+k-a}}\Bigr\}\,.
\end{eqnarray*}
The above relation leads immediately to
\begin{equation}
u_k(0) = k!(k+1)_q\Bigl\{\sum_{i=1}^m(q^{-1}\mu_i)^{k+1} -
\sum_{j=1}^n(q\nu_j)^{k+1}\Bigr\}:=k!(k+1)_q\pi_{k+1}(\mu,\nu)\quad\forall\,
k\ge 0.
\label{u-pi}
\end{equation}
Now, differentiating the relation (\ref{f-prime}) at $y=0$ gives rise to
$$
f_k = (q-q^{-1})\sum_{r=0}^{k-1}{r \choose  k-1} u_r(0)f_{k-r-1}\,,\quad
k\ge1\,,
$$
where we use an obvious condition $f_0 = 1$. On taking into account
relations (\ref{p-f}) and (\ref{u-pi}), we can easily prove the assertion
of the
lemma. Indeed,
\begin{eqnarray*}
kp_k(\mu,\nu)  &=& \frac{f_k}{(q-q^{-1})(k-1)!} \, =\, \frac{1}{(k-1)!}
\,\sum_{r=0}^{k-1}{r \choose  k-1} u_r(0)f_{k-r-1}
\\
&=&\sum_{r=0}^{k-1}\frac{(r+1)_q}{(k-r-1)!}\,\pi_{r+1}(\mu,\nu)\,f_{k-r-1}
\,=\,\sum_{r=1}^{k}r_q\,\pi_r(\mu,\nu)\frac{f_{k-r}}{(k-r)!}
\\
&=& k_q\pi_k(\mu,\nu)+(q-q^{-1})\sum_{r=1}^{k-1}
r_q\,\pi_r(\mu,\nu)\,p_{k-r}(\mu,\nu)\,.
\end{eqnarray*}
At last, the equivalence of (\ref{p-pi}) and (\ref{P-Pi})
is verified by a direct calculation.\hfill\rule{6.5pt}{6.5pt}

\smallskip
Now we are ready to prove the proposition.
We shall verify that relations   (\ref{q-anti}) and (\ref{q-simm})
turn into identities if we substitute $a_r(M)$,
$s_r(M)$ and $p_{k-r}(M)$ by their spectral parameterizations
(\ref{sk-a}), (\ref{sk-s}) and (\ref{pk})--(\ref{dd-i}) and take into account
the equalities (\ref{a-pi}), (\ref{s-pi})  together with (\ref{p-pi}).
For the definiteness we shall verify the relation (\ref{q-anti}).

To simplify the calculation we shall work with the generating
functions: $A(q^{-1}\mu,q\nu|t)$ for $\{a_k(\mu,\nu)\}_{k\geq 0}$,
$\Pi(q^{-1}\mu,q\nu|t)$ for $\{\pi_k(\mu,\nu\}_{k\geq 1}$, and
$P(\mu,\nu|t)$ for $\{p_k(\mu,\nu)\}_{k\geq 1}$.
Substituting the expression of $\Pi(\dots)$ through $A(\dots)$ (see
(\ref{Pi-AS})) into the equality (\ref{P-Pi}) and gathering together
similar terms we get
$$
{d\over dt}\log \Bigl(P(\mu,\nu|t)\,A(q^{-1}\mu,q\nu|-qt)\Bigr)\, =\,
{d\over dt}\log \Bigl(A(q^{-1}\mu,q\nu|-q^{-1}t)\Bigr)\,.
$$
Integrating this equation leads to
$$
P(\mu,\nu|t)\,A(q^{-1}\mu,q\nu|-qt)\, =\, C A(q^{-1}\mu,q\nu |-q^{-1}t)\,.
$$
The obvious boundary condition $P(\mu,\nu|0) = A(q^{-1}\mu,q\nu|0) = 1$
fixes the integration constant $C$ to unity and we come
to a desired relation: the first one from $(\ref{newton-2})$.\hfill
\rule{6.5pt}{6.5pt}

\newcommand\arxiv[1]{
\href{http://arxiv.org/abs/#1}{\tt arXiv:#1}}

\end{document}